\documentclass[11pt,english,a4paper]{amsart}

\usepackage[utf8]{inputenc}
\usepackage[T1]{fontenc}
\usepackage{amsmath,amssymb,amsthm,mathabx,esint,bbm,dsfont,tikz-cd,mathtools}

\usepackage{enumitem}   
\usepackage{amscd}
\usepackage{amsfonts}
\usepackage{mathrsfs}
\usepackage{setspace}
\usepackage{version}
\usepackage{soul,tikz}

\usepackage[english]{babel}
\usepackage[all]{xy}
\usepackage{xcolor}
\usepackage{pict2e}
\usepackage{hyperref}
\usepackage{mathrsfs}

\theoremstyle{plain}
\newtheorem{Th}{Theorem}[section]
\newtheorem{Lem}[Th]{Lemma}
\newtheorem{Cor}[Th]{Corollary}

\newtheorem{Def}[Th]{Definition}

\newtheorem{Rem}[Th]{Remark}

\newtheorem{theorem}{Theorem}[section]

\theoremstyle{definition}

\theoremstyle{remark}

\numberwithin{equation}{section}

\newcommand{\NN}{\mathbb{N}}
\newcommand{\QQ}{\mathbb{Q}}
\newcommand{\KK}{\mathbb{K}}

\newcommand{\ie}{\textit{i.e.},\;}

\newcommand{\SL}{\mathrm{SL}}

\newcommand{\ZZ}{\mathbb{Z}}

\newcommand{\PP}{\mathbb{P}}
\newcommand{\EE}{\mathbb{E}}

\newcommand{\RR}{\mathbb{R}}
\newcommand{\eps}{\varepsilon}
\newcommand{\law}{\mathrm{Law}}
\newcommand\E{\mathbb{E}}
\newcommand\Z{\mathbb{Z}}
\newcommand\R{\mathbb{R}}

\newcommand\N{\mathbb{N}}

\newcommand\GL{\operatorname{GL}}

\newcommand\Mat{\operatorname{Mat}}

\newcommand\Endo{\operatorname{End}}

\newcommand\diag{\operatorname{diag}}
\newcommand\Q{\mathbb{Q}}

\subjclass[2020]{Primary 60B10; Secondary 22F30}

\title{Transient random walks on the space of lattices}

\author{Axel P\'eneau}
\address{
Institut Denis Poisson, 
Université de Tours,
Faculté des Sciences et Techniques
Bâtiment E2,
Parc de Grandmont,
37200 Tours FRANCE}
\email{axel.peneau@univ-tours.fr}

\author{Cagri Sert}
\address{Mathematics Institute, Zeeman Building University of Warwick, Coventry CV4 7AL, U.K.}
\email{cagri.sert@warwick.ac.uk}

\begin{document}

\begin{abstract}
Given $d \geq 2$, we construct a Zariski-dense random walk on the space of lattices $\SL_d(\R)/\SL_d(\Z)$ that exhibits escape of mass. This negates the suggestion of recurrence made by Benoist \cite{Benoist2014RecurrenceOT} (ICM 2014) and  by Bénard-de Saxcé \cite{Benard2022} (also asked in \cite{benoist-quint.tagaki}). For any $p \in (0,1)$, we also construct such a random walk with finite $L^p$-moment which shows that the moment assumption in \cite{Benard2022} is sharp. 
\end{abstract}


\maketitle

\section{Introduction}

Let $G$ be a Lie group and $\Lambda<G$ a lattice, i.e.\ a discrete subgroup of $G$ with finite covolume. A probability measure $\mu$ on $G$ defines a Markov chain on the state space $X=G/\Lambda$: for $x \in X$, the transition probability is given by $\mu \ast \delta_x$ which is the push-forward of $\mu$ under the map $g \mapsto gx$. Similarly, the $n^{th}$-step distribution of this Markov chain starting from $x \in X$ is $\mu^{\ast n} \ast \delta_x$, where $\mu^{\ast n}$ denotes the $n$-fold convolution of $\mu$. We will sometimes refer to this Markov chain as a \textit{random walk on} $X$. When $G<\GL_d(\R)$ is a linear Lie group, we will say that a probability measure $\mu$ (equivalently the induced random walk on $X$) is \textit{Zariski-dense} if the semi-group $\Gamma_\mu$ generated by the support of $\mu$ is Zariski-dense in $G$. 
The random walk on $X$ induced by a probability measure $\mu$ on $G$ is said to be \textit{recurrent in law on} $X$ if for every $\varepsilon>0$ and $x \in X$ there exists a compact $M \subseteq X$ and $n_0 \in \mathbb{N}$ such that $\mu^{\ast n}\ast \delta_x(M)>1-\varepsilon$ for every $n \geq n_0$.


It was suggested in \cite{Benoist2014RecurrenceOT} and in \cite[\S 7]{Benard2022} that when $\Lambda$ is non-uniform, any Zariski-dense random walk is recurrent in law on $X$ (also asked in \cite[Question 2]{benoist-quint.tagaki}). The following result says in particular that this is not the case:

\begin{theorem}[Full escape]\label{thm.intro.full.escape}
Given $d \geq 2$, there exist a Zariski-dense probability measure $\mu$ on $\SL_d(\R)$ and $x \in X=\SL_d(\R)/\SL_d(\Z)$ such that $\mu^{\otimes \N}$-a.s.\  $\gamma_0\cdots \gamma_{n-1}\cdot x \to \infty$ and $\gamma_{n-1}\cdots \gamma_{0}\cdot x \to \infty$ in $X$ as $n \to \infty$. In particular, for every compact subset $M$ of $X$, we have $\mu^{\ast n} \ast \delta_x(M) \to 0$ and hence the $\mu$-random walk on $X$ is not recurrent in law.
\end{theorem}

The homogeneous space $\SL_d(\R)/\SL_d(\Z)$ is also known as the space of lattices in $\R^d$. Indeed, the map $g \mapsto g \cdot \Z^d$ yields a natural identification between $\SL_d(\R)/\SL_d(\Z)$ and the set of lattices in $\R^d$ with unit covolume. With the quotient topology, this space is non-compact (equivalently, the lattice $\SL_d(\Z)$ is non-uniform). The escape to infinity is captured by Mahler's compactness criterion which says that a sequence $x_n \in \SL_d(\R)/\SL_d(\Z)$ escapes to infinity as $n \to \infty$ if and only if the sizes of smallest non-zero vectors in the associated sequence of lattices tend to $0$.

Random walks on homogeneous spaces $G/\Lambda$ in this generality were first studied by Eskin--Margulis \cite{eskin-margulis}. Their result shows in the setting of the theorem above that if $\mu$ has a finite exponential moment, then the random walk is recurrent on $X$ (even uniform in the sense that the compact $M$ can be chosen independently of $x \in X$). Here finite exponential moment means that there exists $\alpha>0$ such that $\int \|g\|^\alpha d\mu(g)<\infty$, where $\|.\|$ is any choice of operator norm on $\Endo(\R^d)$. In their work, Eskin--Margulis employed Foster-Lyapunov recurrence technique tracing back in this context to the work of Eskin--Margulis--Mozes \cite{eskin-margulis-mozes}. This technique was recently improved in this context in the aforementioned work \cite{Benard2022} of Bénard--de Saxcé who managed to prove the recurrence in law of the random walk on $X$ under finite $L^1$-moment condition on $\mu$, i.e.\ $\int \log \|g\| d\mu(g)<\infty$.

Our next result says that, for $p \in (0,\infty]$, among the $L^p$-moment assumptions (i.e.\ $\int (\log \|g\|)^p d\mu(g)<\infty$), the $L^1$-assumption of \cite{Benard2022} is sharp for recurrence in law:

\begin{theorem}[Escape of mass under $L^p$-moment, $p \in (0,1)$]\label{thm.intro.p.moment}
    Given $d \geq 2$ and $p \in (0,1)$, there exists a Zariski-dense probability measure $\mu$ on $\SL_d(\R)$ with finite $L^p$-moment such that the associated random walk on $\SL_d(\R)/\SL_d(\Z)$ is not recurrent in law on $X$. 
\end{theorem}

As opposed to Theorem \ref{thm.intro.full.escape}, in Theorem \ref{thm.intro.p.moment} we do not prove that there is a full escape of mass starting from some $x \in X$ (i.e.\ $\mu^{\ast n} \ast \delta_x(M) \to 0$ as $n \to \infty$ for every compact $M$). We do not know whether there are examples with finite $L^p$-moment assumption exhibiting full escape of mass.

\bigskip


There are two special aspects of the constructions in Theorems \ref{thm.intro.full.escape} and \ref{thm.intro.p.moment}. On the one hand the probability measure $\mu$ driving the random walk is a discrete measure on $\SL_d(\R)$ supported on commensurator $\SL_d(\Q)$ of the lattice $\SL_d(\Z)$. This arithmetic aspect plays a seemingly important role in the proofs. 
On the other hand, and relatedly, the arguments used to prove the above results allow us to guarantee the escape of mass of the random walk starting only from countably many points in $\SL_d(\R) / \SL_d(\Z)$. Although this is enough to 
negate the recurrence in law (suggested in \cite{Benoist2014RecurrenceOT,Benard2022}), one may argue that these points have ``obvious reasons to escape'' in the sense of \cite[\S 1. B]{barak.gafa.divergent} or \cite{dani.divergent}, and the existence of such examples relies on an arithmetic feature. For this reason one is led to ask whether ``non-obvious divergent trajectories'' exist. 

The following two results address both of these aspects demonstrating that they are not essential for the non-recurrence, or in other words, non-obvious divergent trajectories exist. Moreover, their proofs indicate that constructions of non-recurrent random walks do not rely on an arithmetic feature but instead exhibits a certain flexibility with a Diophantine flavor.

\begin{theorem}[Large support non-recurrent random walks]\label{th:non-alg}
    Let $0 < p < 1$ and $d \ge 2$.
    There exist a constant $\eps > 0$ and a symmetric probability distribution $\mu$ with finite $L^p$-moment, having no atoms and with support equal to $\SL_d(\RR)$ such that for all compact $K \subset \SL_{d}(\RR) / \SL_{d}(\ZZ)$, we have
    \begin{equation}\label{eq:lim-transcendental}
        \liminf_{n \to +\infty} \frac{1}{n}\sum_{k = 1}^{n}\mu^{*n}(K \cdot \SL_{d}(\ZZ)) \le 1 - \eps.
    \end{equation}
\end{theorem}

In both above and below results, the symmetry (i.e.\ invariance under push-forward by the map $\gamma \mapsto \gamma^{-1}$) is an additional feature. With an easy modification, one can construct non-symmetric probability measures with the same properties.

\begin{theorem}[Uncountably many divergent points]\label{th:non-alg-escape}
    Let $d \ge 2$. There exists an uncountable, dense and measurable subset $S \subset X$ and a symmetric probability distribution $\mu$ on $\SL_d(\RR)$ which has no atoms, has full support and such that given $(\gamma_n)\sim \mu^{\otimes\NN}$, we have:
    \begin{equation*}
        \PP\left(\forall x \in S,\, \gamma_1\cdots\gamma_n \cdot x \underset{n \to \infty}{\to} \infty \right) = 1.
    \end{equation*}
\end{theorem}

Note that by Chacon--Ornstein Theorem, the set of points $x \in X$ for which the Ces\`{a}ro averages $\frac{1}{n}\sum_{k=1}^n \mu^{\ast k} \ast \delta_x$ exhibit escape of mass has zero measure for the Haar measure on $X$. Due to the specific nature of our constructions, it is possible that in our examples the Hausdorff dimension of divergent points can be analysed using techniques in the study of divergent trajectories under diagonalizable flows (to mention  a few \cite{dani.divergent, cheung.annals, kadyrov-etal, liao-etal}). However, in general we believe that the set of divergent starting points reflects a subtle Diophantine property of the heavy-tails of the probability measure $\mu$. Indeed, we believe that the recurrence suggestion by Benoist \cite{Benoist2014RecurrenceOT} may be true when amended with a subtle Diophantine condition.

\bigskip

Finally, in a different direction, when $\mu$ is not supposed to be Zariski-dense but instead supported on a unipotent subgroup, examples of divergent random walks were found by Breuillard \cite{breuillard.escape} even with compact support (which is impossible under Zariski-density assumption by Eskin--Margulis \cite{eskin-margulis}). The example of $\mu$ in \cite{breuillard.escape} is necessarily non-symmetric. As mentioned above, this aspect is not relevant for us: we have both symmetric and non-symmetric examples for all studied phenomena.

\section{Preliminaries on rational lattices}

The goal of this brief section is to collect some basic facts about the $\SL_d(\Q)$-orbit of $\Z^d$ in the space of lattices and record the key Lemma \ref{lem:systole}.

In what follows, we endow $\RR^d$ with the Euclidean norm and $\SL_d(\RR)$ with the induced operator norm $\|g\| := \max_{x\in\RR^d\setminus\{0\}} \frac{\|g x\|}{\|x\|}$. For $i=1,\ldots,d$, we denote by $e_i$ the $i^{th}$ standard basis vector in $\R^d$.

\begin{Def}[Systole of a lattice]
	Let $\Lambda \subset \RR^d$ be a lattice. We call the quantity $\min\{\|v\| : v \in \Lambda\setminus\{0\}\}$ the systole of $\Lambda$ and denote it by $\delta(\Lambda)$.
\end{Def}

The following elementary linear algebraic fact will be used to show that once the systole is made very small it cannot be enlarged by a relatively medium sized matrix.

\begin{Lem}\label{lem:expand}
	Let $g \in\SL_d(\RR)$ and let $\Lambda \subset \RR^d$ be a lattice. Then $\|g^{-1}\|^{-1} \delta(\Lambda) \le \delta(g\cdot\Lambda) \le \|g\| \delta(\Lambda)$. 
\end{Lem}

\begin{proof}
	Let $x \in\Lambda\setminus\{0\}$ be such that $\|x\| = \delta(\Lambda)$. Then $gx \in g \cdot\Lambda\setminus\{0\}$ so $\delta(g\cdot\Lambda) \le \|gx\| \le \|g\| \delta(\Lambda)$. By the same reasoning applied to the matrix $g^{-1}$ and the lattice $g \cdot \Lambda$, we get that $\delta(\Lambda) \le \|g^{-1}\| \delta(g\cdot\Lambda)$.
\end{proof}

\begin{Def}[Rational height of a matrix]
	Let $g \in \mathrm{SL}_d(\QQ)$. Let $q(g)$ be the smallest common denominator of the entries of $g$. For simplicity, we will also use the notation $q'(g):=q(g^{-1})$. 
\end{Def}

\begin{Lem}\label{lem:inter-lattice}
    Let $g \in \mathrm{SL}_d(\QQ)$ and $x \in \ZZ^d$. Then $q'(g) x \in g \cdot \ZZ^d$. In particular, $g\cdot \Z^d$ contains $r e_i$ with $ 0 \neq |r| \leq q'(g)$ for any $i=1,\ldots,d$.
\end{Lem}
\begin{proof}
   Indeed, $q'(g) x \in g \cdot \ZZ^d$ if and only if $g^{-1}q'(g) x \in \Z^d$. But since $g^{-1}q'(g) \in \Mat_d(\Z)$, the claim follows.
\end{proof}




	

\begin{Def}[Height of a rational matrix]\label{def.height}
	Let $g \in \mathrm{SL}_d(\QQ)$. We define the height of $g$ as:
	\begin{equation*}
		H(g) := \log\max\left\{\|g\|,\|g^{-1}\|, q(g), q'(g)\right\}.
	\end{equation*}
\end{Def}
We remark in passing that $\|g^{-1}\| \leq \|g\|^{d-1}$ and $q'(g) \leq q(g)^{d-1}$ and in particular $H(g) \leq (d-1) \log \max \{\|g\|, q(g)\}$. Moreover, $H$ is subbadditive, i.e.\ for every $g,h \in \SL_d(\Q)$, $H(gh) \leq H(g)+H(h)$ and satisfies $H(g)=H(g^{-1})$.

We now prove the lemma that will allow us to construct transient random walks.

\begin{Lem}\label{lem:systole}
	Let $0 \le k < n$ be integers. 
	Let $(\gamma_0,\cdots,\gamma_{n-1}) \in \mathrm{SL}_d(\QQ)^n$ and let $g := \gamma_0 \cdots \gamma_{n-1}$. 
	Assume that $\gamma_k$ is a diagonal matrix $\diag(2^{-m}, 2^m, 1, \dots, 1)$ for some integer $m \ge 0$. 
	Then:
	\begin{equation}\label{eq:dom-height}
		-\log \delta(g \cdot \ZZ^d) \ge H(\gamma_k) - \sum_{i \neq k} H(\gamma_i).
	\end{equation}
\end{Lem}

\begin{proof}
	We prove that:
	\begin{equation}
		\delta(g \cdot \ZZ^d) \le 2^{-m} \prod_{i= 0}^{k-1} \|\gamma_i\| \prod_{i = k+1}^{n-1} q'(\gamma_i),
	\end{equation}
	from which \eqref{eq:dom-height} follows directly.
	First we have $q'(\gamma_{k+1} \cdots \gamma_{n-1}) \le \prod_{i = k+1}^{n-1} q'(\gamma_i)$ so by Lemma \ref{lem:inter-lattice}, there is a vector $x$ in the segment $[0,  \prod_{i = k+1}^{n-1} q'(\gamma_i) e_1]$ that is in $\gamma_{k+1}\cdots \gamma_n \cdot \ZZ^d$.
	Therefore, there is a vector in the segment $[0, 2^{-m} \prod_{i = k+1}^{n-1} q'(\gamma_i) e_1]$ that is in the lattice $\gamma_{k}\cdots \gamma_n \cdot \ZZ^d$.
	Hence,
	\begin{equation*}
		\delta(\gamma_{k}\cdots \gamma_n \cdot \ZZ^d) \le 2^{-m} \prod_{i = k+1}^{n-1} q'(\gamma_i).
	\end{equation*}
	We then conclude using Lemma \ref{lem:expand}.
\end{proof}

\section{Construction of the rational counterexamples}

 
\subsection{Laws with heavy records}\label{subsec.heavy.record}
In this part, after defining the notion of a probability measure with heavy records, we show the existence of such measures (Lemma \ref{lem:ex-heavy-record}) and a stability property (Lemma \ref{lem:stability}). We then prove Theorem \ref{thm.intro.full.escape} from the introduction.

\begin{Def}[Heavy records]
	Let $\eta$ be a non-trivial probability distribution on $\RR_{\ge 0}$. Let $(x_n)_{n\in\NN} \sim \eta^{\otimes \NN}$. We say that $\eta$ has heavy records if we have almost surely:
	\begin{equation*}
		\frac{\max_{i = 1}^{n} x_i}{\sum_{i = 1}^{n} x_i} \underset{n \to + \infty}{\longrightarrow} 1.
	\end{equation*}
\end{Def}

\begin{Lem}[Stability under small perturbation] \label{lem:stability}
	Let $X$ and $Y$ be real random variables defined on the same probability space. 
	Assume that $Y$ is in $\mathrm{L}^1$ and that both $X$ and $X + Y$ are non-negative.
	Then $\law(X)$ has heavy records if and only if $\law(X + Y)$ has.
\end{Lem}

\begin{proof}
	Note that $X$ and $X + Y$ play symmetric roles so we only need to show one implication. 
	Assume that the law of $X$ has heavy records. 
	Let $(x_n,y_n)_{n\in\NN} \sim \law(X, Y)^{\otimes \NN}$.
Observe first that if a random variable $Z$ belongs to $L^1$ and $(z_n)_{n \in \N} \sim \law(Z)^{\otimes \N}$, then by the law of large numbers $\max_{i=1}^n z_i=o(n)$. It follows that if for a random variable $Z$, $\law(Z)$ has heavy records, then $\mathbb{E}[Z]=+\infty$. In particular, $E[X]=+\infty$, Moreover, since $Y$ is in $L^1$ and $\law(X)$ has heavy records, $\law(X+Y)$ is non-trivial. 
Now, the fact that $\E[X]=+\infty$ together with the law of large numbers imply that $n=o(\sum_{i=1}^n x_i)$ and  hence $\sum_{i=1}^n y_i=o(\sum_{i=1}^n x_i)$ a.s. Moreover, since $Y \in L^1$, $\max_{i=1}^n y_i=o(n)$ and by the heavy records property $n=o(\max_{i=1}^n x_i)$ a.s. The result follows. 	
\end{proof}

Two direct consequence of Lemma \ref{lem:stability} is that if a random variable $X$ has a law with heavy records, then $\EE(X) = +\infty$ and the law of the integer part $\lfloor X \rfloor$ also has.

We now construct a law with heavy records.

\begin{Lem}\label{lem:ex-heavy-record}
	Let $\eta$ be the push-forward of the uniform probability distribution on $[0, 1]$ by the map $t \mapsto \exp(t^{-2})$. Then $\eta$ has heavy records.
\end{Lem}

\begin{proof}
	Let $(t_n)_{n \ge 1}$ be a sequence of independent uniformly distributed random variables with values in $[0,1]$. 
	For all $n \in \NN$ and for all $1 \le k \le n$, let $u_{k,n}$ be the $k$-th smallest value of $\{t_1,\dots, t_n\}$.
	We want to show that almost surely, we have:
	\begin{equation}
		\frac{\max_{i = 1}^{n} \exp(t_i^{-2})}{\sum_{i = 1}^{n} \exp(t_i^{-2})} = \frac{\exp(u_{1,n}^{-2})}{\sum_{k = 1}^{n} \exp(u_{k,n}^{-2})} \underset{n \to + \infty}{\longrightarrow} 1.
	\end{equation}
	For that, by Borel--Cantelli Lemma, it is enough to show that:
	\begin{equation}
		\forall \eps > 0,\; \sum_{n = 1}^{+ \infty}\PP\left({\eps \exp(u_{1,n}^{-2})}\le {n\exp(1/u_{2,n}^{-2})}\right) < + \infty.
	\end{equation}
	Let $\eps > 0$ and $n\in\NN$ be fixed. 
	We have ${\eps \exp(u_{1,n}^{-2})}\le {n\exp(u_{2,n}^{-2})}$ if and only if $u_{1,n}^{-2} - u_{2,n}^{-2} \le \log(n/\eps)$.
	First note that if $u_{1,n} \le n^{-3/4}$ and $u_{2,n} \ge {u_{1,n}} + n^{-9/4}$, then for all $n \ge 2$, we have:
	\begin{equation*}
		 u_{1,n}^{-2} - u_{2,n}^{-2} \ge n^{6/4}\left(1 - \frac{1}{(1+n^{-5/4})^2}\right)
		 \ge n^{1/4}.
	\end{equation*}
	Indeed $1-\frac{1}{(1+x)^2} = \frac{2x + x^2}{1 + 2x + x^2} \ge x$ for all $x \le \frac{1}{2}$. 
	Moreover, there is an integer $k_\eps$ such that for all $n \ge k_\eps$, we have $n^{1/4}\ge \log(n/\eps)$.
	Therefore, we only need to show that $\sum \PP\left(u_{1,n} \ge n^{-3/4}\right) < + \infty$ and $\sum \PP\left(u_{2,n} \le {u_{1,n}} + n^{-9/4}\right) < + \infty$. The first one is straightforward: for all $n \ge 1$, we have
	\begin{equation*}
		\PP\left(u_{1,n} \ge n^{-3/4}\right) = \PP\left(\forall j \le n,\,t_j \ge n^{-3/4} \right) = \left(1 - n^{-3/4}\right)^n \le \exp(-n^{1/4}).
	\end{equation*}
	It therefore remains to show that $\sum \PP\left(u_{2,n} \le {u_{1,n}} + n^{-9/4}\right) < + \infty$. Note that for all $i \neq j$, the conditional distribution of $t_j$ relatively to $t_i$ knowing that $t_i < t_j$ is uniform in $[t_i, 1]$. Therefore, for all $n \in\NN$, the distribution of $u_{2,n}$ relatively to $u_{1,n}$ is the minimum of $(n-1)$ random variables that are independent and uniformly distributed in $[u_{1,n}, 1]$. Therefore we have:
    \begin{equation*}
        \PP\left(u_{2,n} \le {u_{1,n}} + n^{- 9/4}\,\middle|\,u_{1,n}\right) \le (n-1) \left(\frac{n^{-9/4}}{1-u_{1,n}}\right).
    \end{equation*}
    Now we differentiate the cases $u_{1,n} \le \frac{1}{2}$ and $u_{1,n} \ge \frac{1}{2}$:
    \begin{align*}
        \PP\left(u_{2,n} \le {u_{1,n}} + n^{-9/4}\right) 
        & \le \PP\left(u_{2,n} \le {u_{1,n}} + n^{-9/4}\,\middle|\,u_{1,n} \le \frac{1}{2}\right) + \PP\left(u_{1,n} \ge \frac{1}{2}\right) \\
        & \le 2 n^{-5/4} + 2^{-n}.
    \end{align*}
    This concludes the proof.
\end{proof}



\begin{proof}[Proof of Theorem \ref{thm.intro.full.escape}]
	Let $\kappa$ be a finitely supported probability measure on $\SL_d(\QQ)$ such that $\Gamma_\kappa$ is Zariski-dense in $\SL_d(\R)$. 
	Let $\eta$ be the push-forward of the uniform probability distribution on $[0, 1]$ by the map $t \mapsto \lfloor\exp(t^{-2})\rfloor$. 
	Then $\eta$ has heavy records by Lemma \ref{lem:stability} and Lemma \ref{lem:ex-heavy-record}.
	Let $\nu$ be the push-forward of $\eta$ by the map:
	\begin{equation}
		n \longmapsto \mathrm{diag}(2^n,2^{-n},1,\dots,1).
	\end{equation}
	Let $\mu = (\kappa + \nu) / 2$. 
	Then $H_*\nu = \eta$ has heavy records and $H_*\kappa$ is bounded, where $H$ is the height function (Definition \ref{def.height}).
	Therefore, $H_*\mu$ has heavy records by Lemma \ref{lem:stability}.
	
	Now let $(\gamma_n)\sim \mu^{\otimes\NN}$.
	For all $n \in\NN$, write $m_n = \max_{0 \le k < n} H(\gamma_k)$ an $s_n = \sum_{k = 0}^{n-1} H(\gamma_k)$. 
	Then $-\log\delta(\gamma_0\cdots \gamma_{n-1}\cdot\ZZ^d) \ge 2m_n - s_n$ by Lemma \ref{lem:systole}. 
	As a consequence $\delta(\gamma_0\cdots \gamma_{n-1}\cdot\ZZ^d) \to 0$ almost surely.
	With the same argument, we also show that $\delta(\gamma_{n-1}\cdots \gamma_{0}\cdot\ZZ^d) \to 0$ almost surely.
\end{proof}

\subsection{Laws with escape of mass and finite polynomial moment}\label{subsec.Lp}

Here, we introduce the notion of a probability measure with escape of mass, a weakening of heavy records property. We prove the existence of such measures (Lemma \ref{lem:ex-lp}) and after recording a stability property (Lemma \ref{lem:compo-lp}), we prove Theorem \ref{thm.intro.p.moment}.

\begin{Def}
    Let $\eta$ be a probability distribution on $\RR_{\ge 0}$. 
    Let $(x_n)_{n\in\NN} \sim \eta^{\otimes \NN}$. We say that $\eta$ has escape of mass if:
    \begin{equation*}
        \exists \eps >0, \; \forall M \in\RR, \; \limsup_{n \to + \infty}\PP\left(2\max_{1 \le i \le n} x_i - \sum_{i = 1}^n x_i \ge M\right) \ge \eps.
    \end{equation*}
\end{Def}

Note that a probability measure with heavy records has escape of mass and one with finite $L^1$-moment has no escape of mass by the law of large numbers. 

In the next lemma, for every $p' \in (0,1)$, we construct a law with finite $L^{p'}$-moment and escape of mass.

\begin{Lem}\label{lem:ex-lp}
    Let $p > 1$ and let $\eta$ be the push-forward of the uniform probability distribution on $[0,1]$ by the map $t \mapsto t^{-p}$. Then $\eta$ has escape of mass. More precisely, there exists a constant $\eps_p > 0$ such that given $(t_n)_{n\ge 1}$ independent and uniformly distributed in $[0,1]$, we have:
    \begin{equation}\label{eq:eps-p}
        \forall n \ge 1,\,\PP\left(2\max_{1 \le i \le n} t_i^{-p} - \sum_{i = 1}^n t_i^{-p} \ge (2n)^p\right) \ge \eps_p.
    \end{equation}
\end{Lem}

\begin{proof}
    Let $(t_n)_{n\ge 1}$ be independent and uniformly distributed random variables on $[0,1]$.
    For all $1 \le k \le n$, let $u_{k,n}$ be the $k$-th smallest element of $\{t_1, \dots, t_n\}$.
    We want to exhibit a constant $\eps_p > 0$ such that \eqref{eq:eps-p} holds.
    First note that for all $n \in\NN$, we have:
    \begin{equation*}
        2\max_{1 \le i \le n} t_i^{-p} - \sum_{i = 1}^n t_i^{-p} = u_{1,n}^{-p} - \sum_{k = 2}^n u_{k,n}^{-p}.
    \end{equation*}
    The strategy of proof goes as follows. First we show that there exists a universal constant $\alpha >0$ such that:
    \begin{equation}\label{eq:lln}
        \forall n \ge 1, \;  \; \PP\left(\forall \, 2 \le k \le n,\; u_{k,n} \ge \frac{k}{2n}\right) \ge \alpha.
    \end{equation}
    Then we show that there is a constant $\eps_p$ such that:
    \begin{equation}\label{eq:eps-ldev}
        \forall n \ge 1, \; \; \PP\left(\forall \, 2 \le k \le n,\; \frac{k}{2n} \le u_{k,n} \; \; \text{and} \; \; u_{1,n} \le \frac{1}{2n}\left(\sum_{\ell = 1}^{\infty} \ell^{-p}\right)^{-\frac{1}{p}}\right) \ge \eps_p.
    \end{equation}
    Note that if $u_{k,n}\le \frac{k}{2n}$ for all $2\le k \le n$, then $\sum_{k = 2}^n u_{k,n}^{-p} \le (2n)^p\sum_{k = 2}^{\infty} k^{-p}$. Note also that if $u_{1,n} \le \frac{1}{2n}\left(\sum_{k = 1}^{\infty} k^{-p}\right)^{-\frac{1}{p}}$, then $u_{1,n}^{-p} \ge (2n)^p\sum_{k = 2}^{\infty} k^{-p} + (2n)^p$. 
    Therefore \eqref{eq:eps-ldev} implies \eqref{eq:eps-p}.
    First we prove \eqref{eq:lln}. Note that for all $k \le n$, we have:
    \begin{align*}
        \PP\left(u_{k,n} < \frac{k}{2n}\right) & = \PP\left(\# \left\{1\le i \le n\,\middle|\,t_i\le \frac{k}{2n}\right\} \ge k\right)  \\
        & \le \PP\left(\prod_{i = 1}^n \exp(\mathds{1}_{t_i \le \frac{k}{2n}}/10)> \exp(k/10)\right)\\
        & \le \exp(-k/10) \EE\left(\exp\left(\mathds{1}_{t_1 \le \frac{k}{2n}}/10\right)\right)^n\\
        & \le \exp(-k/10) \left(1 + \frac{k}{2n}(\exp(1/10) - 1)\right)^n \\
        & \le \exp\left(-k/10\right)\exp\left(\frac{\exp(1/10) - 1}{2}\right) \le \exp(-k/30),
    \end{align*}
    where we used the i.i.d. property of $t_i$'s in the second inequality and $(1+\alpha)^n \leq e^{\alpha n}$ in the second to last inequality.
    Note moreover that for all $k \le n$, one has:
    \begin{equation}
        \PP\left(u_{k,n} \ge \frac{k}{2n}\,\middle|\, \forall \, 2\le j<k,\,u_{j,n}\ge\frac{j}{2n}\right) \ge \PP\left(u_{k,n} \ge \frac{k}{2n}\right) \ge 1-\exp(-k/30).
    \end{equation}  
    Then by induction, we show that:
    \begin{multline}
        \PP\left(\forall \, 2\le j \le n,\,u_{j,n}\ge\frac{j}{2n}\right) = \prod_{k = 2}^{n}\PP\left(u_{k,n} \ge \frac{k}{2n}\,\middle|\, \forall \, 2\le j<k,\,u_{j,n}\ge\frac{j}{2n}\right) \\ \ge \prod_{k=2}^{n}(1-\exp(-k/30)).
    \end{multline}
    Therefore, we have \eqref{eq:lln} for $\alpha = \prod_{k = 2}^{\infty}(1-\exp(-k/30)) > 0$.
	
    Write $a_p := \left(\sum_{k = 1}^{\infty} k^{-p}\right)^{-\frac{1}{p}}$. Note that the conditional distribution of $u_{1,n}$ relatively to $(u_{k,n})_{k\ge 2}$ is simply the uniform probability distribution in $[0,u_{2,n}]$, therefore :
    \begin{equation}\label{eq:condition}
        \PP\left(u_{1,n} \le \frac{a_p}{2n}\,\middle|\,(u_{k,n})_{k\ge 2}\right) = \min\{1, \frac{a_p}{2n u_{2,n}}\}.
    \end{equation}
	Note also that for all $K \ge 0$ and for all $n > K$, one has:
	\begin{equation}
		\PP\left(u_{2,n} \ge \frac{K}{n}\right) = \left(1-\frac{K}{n}\right)^n + \binom{n}{1} \frac{K}{n} \left(1-\frac{K}{n}\right)^{n-1} \le \left(1 + \frac{K}{1-K/n}\right) \exp(-K).
	\end{equation}
	Let $K$ be such that $(1+2K) \exp(-K) \le \alpha/2$. Then for all $n \ge 2K$, one has $\PP\left(u_{2,n} \ge \frac{K}{n}\right) \le \frac{\alpha}{2}$ and for all $n \le 2K$, one has $\PP\left(u_{2,n} \ge \frac{2K}{n}\right) = 0$. Therefore, for all $n \in\NN$, one has:
	\begin{equation*}
		\PP\left(u_{2,n} \le \frac{2K}{n}\,\middle|\,\forall \, 2\le j \le n,\,u_{j,n}\ge\frac{j}{2n}\right) \ge \frac{1}{2}.
	\end{equation*}

	Then by \eqref{eq:condition}, one has:
	\begin{equation*}
		\PP\left(u_{1,n} \le \frac{a_p}{2n}\,\middle|\,\forall \, 2\le j \le n,\,u_{j,n}\ge\frac{j}{2n}\right) \ge \frac{a_p}{8K}.
	\end{equation*}
	As a consequence, we have \eqref{eq:eps-ldev} and therefore \eqref{eq:eps-p} for $\eps_p = \frac{\alpha a_p}{8K}$.  
\end{proof}

We remark in passing that for $p>1$, the probability measure $\eta$ considered in the previous lemma has a finite $L^{p'}$-moment if and only if $p'<\frac{1}{p}$, and $a_p$ (and hence $\varepsilon_p$) converges to zero as $p \searrow 1$.

\begin{Rem}
	For a more intuitive understanding of lemma \ref{lem:ex-lp}, note that when we rescale, the family of random variables $(nu_{1,n},\dots, nu_{n,n})$ converges for the weak-$*$ topology to a Poisson clock.
	Let $(t_i)$ be independent and uniformly distributed random variables in $[0,1]$ and let $(y_k)_{k \ge  1}$ be ticking times of the Poisson clock of parameter $1$ in $\RR_{\ge 0}$.
	Then for all $M \in\RR$, one has:
	\begin{equation*}\label{eq:schlob}
		\lim_{n \to + \infty}\PP\left(2\max_{1 \le i \le n} t_i^{-p} - \sum_{i = 1}^n t_i^{-p} \ge M\right) = \PP\left(y_1^{-p} \ge \sum_{k = 2}^{+\infty} y_k^{-p} \right).
	\end{equation*}
	Note that for $p \le 1$, the sum $\sum_{k = 2}^{+\infty} y_k^{-p}$ is almost surely not defined.
\end{Rem}

\begin{Lem}\label{lem:compo-lp}
    Let $p > 1$ and let $\eta$ be the push-forward of the uniform probability distribution on $[0,1]$ by the map $t \mapsto t^{-p}$. 
    Let $\eps_p$ be as in Lemma \ref{lem:ex-lp}. 
    Let $\alpha \in (0,1)$, $M \ge 0$, and $\kappa$ be a probability distribution on $[0, M]$.
    Let $\mu = \alpha \eta + (1-\alpha) \kappa$ and  $(x_n) \sim \mu^{\otimes\NN}$.
    Then we have:
    \begin{equation}\label{eq:eps-p-general}
        \forall n \ge 1,\,\PP\left(2\max_{1 \le i \le n} x_i - \sum_{i = 1}^n x_i \ge \alpha^p n^p - nM\right) \ge \frac{\alpha\eps_p}{2}.
    \end{equation}
\end{Lem}

\begin{proof}
    Let $\mathcal{B}_{\alpha} = \alpha \delta_1 + (1-\alpha) \delta_0$ be the Bernoulli distribution of parameter $\alpha$.  
    Let $(a_n)\sim \eta^{\otimes \NN}$, let $(b_n) \sim \kappa^{\otimes\NN}$ and let $(i_n)\sim \mathcal{B}_{\alpha}^{\otimes\NN}$ be globally independent.
    Then $(i_n a_n + (1- i_n) b_n)_n \sim \mu^{\otimes\NN}$ so we write $x_n := i_n \alpha_n + (1- i_n) b_n$ for all $n \in\NN$. 
    Let $s_n := i_1 + \dots + i_n$ for all $n \in\NN$.
    We claim that $\PP(s_n \ge \alpha n/2) \ge \alpha/2$. 
    Indeed $\EE(n- s_n) = (1-\alpha)n$ so by Markov's inequality, we have:
    \begin{equation*}
        \PP(s_n < \alpha n/2) = \PP\left(\frac{n-s_n}{2} > 1-\alpha/2\right) \le \frac{1-\alpha}{1-\alpha/2} \le 1-\alpha/2.
    \end{equation*}
    For all $k \in \NN$, we write $j_k := \min\{j \ge 1\,|\, s_j = k\}$ and $y_k = x_{j_k}$.
    Then the random integer $j_k$ is almost surely well-defined and $s_{j_k - 1} < s_{j_k}$ so $i_{j_k} = 1$ for all $k$ and therefore $y_k = a_{j_k}$ for all $k$.
    The data of $(j_k)_{k}$ is independent of the data of $(a_j)_j$ so $(y_k) \sim \eta^{\otimes\NN}$ and the joint data of $(y_k)_k$ is independent of the joint data of $(s_n)_n$.
    By Lemma \ref{lem:ex-lp}, we have:
    \begin{equation*}
        \forall n \ge 1,\,\PP\left(2\max_{1 \le k \le s_n} y_k - \sum_{k = 1}^{s_n} y_k \ge (2s_n)^p\,\middle|\, s_n\right) \ge \eps_p.
    \end{equation*}
    Therefore
    \begin{equation*}
        \forall n \ge 1,\,\PP\left(2\max_{1 \le k \le s_n} y_k - \sum_{k = 1}^{s_n} y_k \ge (2s_n)^p \; \; \text{and} \; \; s_n \ge \frac{\alpha n}{2}\right) \ge \frac{\alpha\eps_p}{2},
    \end{equation*}
    which implies \eqref{eq:eps-p-general}.
\end{proof}

We are now in a position to prove Theorem \ref{thm.intro.p.moment} from the introduction.

\begin{proof}[Proof of Theorem \ref{thm.intro.p.moment}]
Let $p \in (0,1)$ and $d \geq 2$ be given. Let $\kappa$ be a probability measure with finite support on $\SL_d(\QQ)$ whose support generates a Zariski-dense subsemigroup in $\SL_d(\R)$. Let $M := \max_{g \in\mathbf{supp}(\kappa)} H(g)$, where we recall that $H(.)$ is the height function from Definition \ref{def.height}.
    Let $p' \in (1,\frac{1}{p})$ and $\eta$ be the push-forward of the uniform probability distribution on $[0, 1]$ by the map $t \mapsto \lfloor t^{-p'}\rfloor$ and $\nu$ be the push-forward of $\eta$ by the map:
    \begin{equation*}
	n \longmapsto \mathrm{diag}(2^n,2^{-n},1,\dots,1).
    \end{equation*}
    Let $\mu = (\kappa + \nu) / 2$. Note that $H_*(\mu)$ is  weakly $\mathrm{L}^{1/p'}$ so it is $\mathrm{L}^{p}$ and in particular $\mu$ has finite $L^p$-moment. Let now $\eps_{p'}>0$ be as in Lemma \ref{lem:ex-lp} and $(\gamma_n)_n \sim \mu^{\otimes\NN}$. 
    Then by Lemma \ref{lem:compo-lp}, applied with $\eta$ and $H_*\kappa$ with $\alpha = 1/2$, we have:
    \begin{equation}
        \forall n \ge 1,\,\PP\left(2\max_{1 \le i \le n} H(\gamma_i) - \sum_{i = 1}^n H(\gamma_i) \ge (n/2)^{p'} - nM\right) \ge \frac{\eps_{p'}}{4}.
    \end{equation}
    So by Lemma \ref{lem:systole}, we have:
    \begin{equation}\label{eq.old.corol.lp}
        \forall n \ge 1,\,\PP\left(- \log\delta(\gamma_1\cdots \gamma_n \cdot \ZZ^d) \ge (n/2)^{p'} - nM\right) \ge \frac{\eps_{p'}}{4}. \qedhere
    \end{equation}
It now follows from Mahler's compactness principle that the random walk induced by $\mu$ on $\SL_d(\R)/\SL_d(\Z)$ is not recurrent in law. 
\end{proof}

\section{Diophantine counter-examples with no atoms}

This final section is devoted to the proofs of Theorems \ref{th:non-alg} and \ref{th:non-alg-escape}.

\subsection{Preliminaries}

First let us quantify how small perturbations influence the systole.

\begin{Lem}\label{lem:perturb-stystole}
    Let $d$ be a positive integer. Let $g , h \in \SL_d(\RR)$.
    Then:
    \begin{equation*}
        \delta((g+h)\cdot\ZZ^d) \le \left({\|h\|}{\|g^{-1}\|} + 1\right) \delta(g \cdot \ZZ^d).
    \end{equation*}
\end{Lem}

\begin{proof}
    Let $x \in \ZZ^d$ be such that $\delta(g\cdot\ZZ^d) = \|g x\|$. 
    Then $\|x\| \le \|g^{-1}\| \|gx\| = \|g^{-1}\| \delta(g\cdot\ZZ^d)$ so $\|hx\| \le \|h\| \|g^{-1}\| \delta(g\cdot\ZZ^d)$. Therefore $\|(g+h) x\| \le \delta(g\cdot\ZZ^d) + \|h\| \|g^{-1}\| \delta(g\cdot\ZZ^d)$
\end{proof}

\begin{Lem}\label{lem:perturb-prod}
    Let $n, d \ge 2$ be integers. 
    Let $\gamma_1, \dots, \gamma_n \in \SL_d(\RR)$.
    Let $\gamma'_1, \dots, \gamma'_n \in \mathrm{Mat}_d(\RR)$ and let $\tau = \max_{1 \le k \le n} \frac{\|\gamma'_k - \gamma_k\|}{\|\gamma_k\|}$. Then:
    \begin{equation*}
        \left\|\prod_{k = 1}^n\gamma'_k - \prod_{k = 1}^n \gamma_k\right\| \le \tau n (1+\tau)^n \prod_{k = 1}^n \|\gamma_k\|.
    \end{equation*}
\end{Lem}

\begin{proof}
    We prove this result by induction on $n$.
    Assume that Lemma \ref{lem:perturb-prod} holds for $n-1$ \ie
    \begin{equation*}
        \left\|\prod_{k = 1}^{n-1}\gamma'_k - \prod_{k = 1}^{n-1} \gamma_k\right\| \le \tau(n-1)(1+\tau)^{n-1} \prod_{k = 1}^{n-1} \|\gamma_k\| .
    \end{equation*}
    Then we multiply everything by $\gamma_n$ and we have:
    \begin{equation*}
        \left\|\gamma_n \prod_{k = 1}^{n-1}\gamma'_k - \prod_{k = 1}^{n} \gamma_k\right\| \le  \tau(n-1) (1+\tau)^{n-1}\prod_{k = 1}^{n}\|\gamma_k\|.
    \end{equation*}
    Moreover, we have:
    \begin{equation*}
        \left\|\prod_{k = 1}^{n}\gamma'_k - \gamma_n \prod_{k = 1}^{n-1}\gamma'_k\right\| \le \|\gamma'_n - \gamma_n\| \prod_{k = 1}^{n-1} \|\gamma'_k\| \le \tau (1+\tau)^n\prod_{k = 1}^{n}\|\gamma_k\|.
    \end{equation*}
    We conclude using the triangular inequality.
\end{proof}

\begin{Lem}\label{lem:perturb-conclusion}
    Let $n, d \ge 2$ be integers. 
    Let $\gamma_1, \dots, \gamma_n \in \SL_d(\RR)$.
    Let $\gamma'_1, \dots, \gamma'_n \in \SL_d(\RR)$ and let $\tau = \max_{1 \le k \le n} \frac{\|\gamma'_k - \gamma_k\|}{\|\gamma_k\|}$.
    Then:
    \begin{equation*}
        \delta\left(\gamma'_1 \cdots \gamma'_n \cdot \ZZ^d\right) \le \left(1 + \tau n (1+\tau)^n \prod_{k = 1}^n \|\gamma_k\|\|\gamma_k^{-1}\|\right) \delta\left(\gamma_1 \cdots \gamma_n \cdot \ZZ^d\right)
    \end{equation*}
\end{Lem}

\begin{proof}
    Let $g := \gamma_1 \cdots \gamma_n$ and let $h := \gamma'_1 \cdots \gamma'_n - \gamma_1 \cdots \gamma_n$.
    By Lemma \ref{lem:perturb-prod}, we have $\|h\| \le \tau n (1+\tau)^n \prod_{k = 1}^n \|\gamma_k\|$.
    Moreover, we have $\|g\|\le \prod_{k = 1}^n \|\gamma_k^{-1}\|$.
    We conclude using Lemma \ref{lem:perturb-stystole}.
\end{proof}

Now we can construct the measure $\mu$ mentioned in Theorem \ref{th:non-alg}.
Let $d \ge 2$ be fixed. 
Let $I_d \in \SL_d(\RR)$ be the identity matrix.
Let $\left(E_{i,j}\right)_{1\le i,j \le d}$ be the canonical basis of $\mathrm{Mat}_{d\times d}(\RR)$.
Given $t \in \RR$, we write $N_t = I_d + t E_{1,2}$.
Note that $t \mapsto N_t$ is a continuous group homomorphism from $\RR$ to $\mathrm{SL}_d(\RR)$.

\subsection{Proof of Theorem \ref{th:non-alg}}

The idea is to construct a measure that has rational approximations that satisfy an escape property similar to \eqref{eq.old.corol.lp} for all $n$ large enough depending on the level of approximation. Around a certain interval of time depending on the level of approximation, the upper bound of the Ces\`{a}ro averages we get using \eqref{eq:eps-p-general} stays below $1- \eps$ before the error term in the rational approximation becomes non-negligible.
The issue is that then we do not have an upper bound on the Ces\`{a}ro averages at all times, so we only get a $\liminf$ in \eqref{eq:lim-transcendental}.



\begin{proof}[Proof of Theorem \ref{th:non-alg}]
    Let $d \ge  2$ be an integer, let $0 < p' < 1$ and let $1 < p < 1/p'$.
    Let $\eta$ be the push-forward of the Lebesgue measure on $[0,1/2]$ by the map $t \mapsto \lfloor t^{-p}\rfloor$.
    Let $\nu$ be the push-forward of $\eta$ by the map $a \mapsto \mathrm{diag}(2^{a}, 2^{-a}, 1, \dots, 1)$.
    It is straightforward to construct a probability measure with full support on $\mathrm{SL}_d(\QQ)$ such that its push-forward by the height function $H$ has finite first moment. So let $\kappa_0$ be such a probability measure. Let $\mu_0 = \nu + \kappa_0/2$ and let $\tilde{\mu}_0 := \frac{\mu_0 + \check{\mu}_0}{2}$, where for a probability measure $\eta$ on $\SL_d(\R)$, we write $\check{\eta}$ for the push-forward of $\eta$ by the inverse map $\gamma \mapsto \gamma^{-1}$.
    The measure $\tilde{\mu}_0$ can alternatively be constructed as follows.
    Let $u \in [0,1]$ be uniformly distributed, let $s \in\{-1,1\}$ be uniformly distributed and independent of $u$ and let $g \sim \kappa_0$ be independent of the joint data of $u$ and $s$. 
    If $u > 1/2$, let $\gamma = g^s$ otherwise, let $\gamma = \mathrm{diag}(2^{s u^{-p}}, 2^{-s u^{-p}}, 1, \dots, 1)$.
    Then $\gamma \sim \tilde{\mu}_0$.
    We write $\mathcal{U}$ for the Lebesgue measure on $[0,1]$.
    
    Let $j > 0$ be an integer and assume that we have constructed a family of exponents $1 \le i_1 < \dots < i_j$.
    Let $\lambda_j$ be the uniform distribution on $\left\{\sum_{k=1}^j\epsilon_{k}2^{-i_k}\,\middle|\, \epsilon_k \in\{0,1\} \right\}$.
    Note that $\lambda_j$ is the push-forward of the uniform probability distribution on $\{0,1\}^\NN$ -- denote it by $\lambda'$ -- by the map $l_j:(\epsilon_k)_{k \ge 1}  \mapsto \sum_{k=1}^j\epsilon_{k}2^{-i_k}$.
    We construct $i_{j+1}$ as follows.
    Let $\mu_j :=  (N_*\lambda_j) * \mu_0$, $\tilde{\mu}_j=\frac{1}{2}(\mu_j + \check{\mu}_j)$, and $\kappa_j := (N_*\lambda_j) * \kappa_0$, where $N$ is the morphism $\R \to \SL_d(\R)$ defined above.
    Let $(u_n, s_n, g_n, k_n)_n \sim \left(\mathcal{U}\otimes(\frac12\delta_{-1}+\frac12\delta_{1})\otimes \kappa_0 \otimes \lambda'\right)^{\otimes\NN}$.
    For all $n \ge 1$, write:
    \begin{equation}\label{eq.def.gamma_n}
        \gamma_n := \mathds{1}_{u_n > 1/2} g_n + \mathds{1}_{u_n \le 1/2} \mathrm{diag}(2^{u_n^{-p}}, 2^{u_n^{-p}}, 1, \dots, 1),
    \end{equation}
    write $h^j_n := N_{l_j(k_n)}$, write $\gamma^j_n = h^j_n \gamma_n$ and write $\tilde{\gamma}^j_n := (\gamma^j_n)^{s_n}$.
    Then $(\gamma_n)\sim\mu_0^{\otimes\NN}$ and $(\gamma^j_n)_n \sim \mu_j^{\otimes\NN}$ and $(\tilde\gamma^j_n)_n \sim \tilde\mu_j^{\otimes\NN}$.
    Note that $H(h^j_n) \le i_j \log(2)$ almost surely and for all $j$ and $n$.
    Let $M := \int H(\gamma) d\kappa_0(\gamma)$ and let $M' := \int \log(\|\gamma\|\|\gamma^{-1}\|)^{p'} d\mu_0(\gamma)$.
    The constant $M$ is finite by construction of $\kappa_0$ and $M'$ is finite by construction of $\nu$.
    Let $a_{j}$ be the smallest integer such that $a_{j}^p \ge (4M / \eps_p + 2i_j \log(2)) a_{j}$, where $\eps_p>0$ is the constant given by Lemma \ref{lem:ex-lp}.
    Let $i_{j+1} \ge i_j + 1$ be the smallest integer such that:
    \begin{equation*}
        2^{-i_{j+1}+1} 2a_{j} (1 + 2^{-i_{j+1}+1})^{2a_j} \exp\left(\left(\frac{4a_{j} (M'+1)}{\eps_p}\right)^{1/p'}\right) \le 1.
    \end{equation*}

    This allows us to define a sequence $(i_j)_{j \in\NN}$ such that $1 \leq  i_1 < i_2 < \dots$. Furthermore, 
    let $l_\infty : (\epsilon_k)_{k \ge 1}  \mapsto \sum_{k=1}^\infty\epsilon_{k}2^{-i_k}$ and $h^\infty_n := N_{l_\infty(k_n)}$, write $\gamma^\infty_n = h^\infty_n \gamma_n$ and write $\tilde{\gamma}^\infty_n := (\gamma^\infty_n)^{s_n}$.
    The distribution of $\gamma^\infty_n$ has no atoms and full support since it is given by the convolution of two probability measures, namely $(N_{\ell_\infty(.)})_\ast \lambda'$ and $\mu_0$, where the former has no atoms and the latter has full support, and each of these properties is preserved under convolution. It is also straightforward to see from the construction that the distribution of $\gamma^\infty_n$ has finite $L^{p'}$-moment. These same properties hold for the distribution of $\tilde{\gamma}^\infty_n$ which is moreover symmetric.  
    
 We now aim to show that for all $j \geq 1$
    \begin{equation}\label{eq.thm.1.3.goal}
        \frac{1}{2a_j}\sum_{n = 1}^{2a_j}\PP\left(-\log\delta\left(\tilde\gamma_1^\infty\cdots \tilde\gamma_n^\infty \cdot \ZZ^d\right) \ge (2^p - 1) a_j^p - \log(2) \right) \ge \frac{\eps_p}{4}.
    \end{equation}
    This easily finishes the proof. Indeed since $i_{j} \to + \infty$ and $a_j \ge (i_j\log(2))^{\frac{1}{p-1}}$, we have $a_j \to + \infty$.
    Therefore, letting $K \subset \SL_d(\RR)/\SL_d(\ZZ)$ be compact, $B := \max_{x \in \KK}(-\log\delta(x))$ and $j_0 := \min\{j \ge 1\,|\,(2^p - 1) a_j^p - \log(2) \ge B\}$,
    for all $j \ge j_0$, we have:
    \begin{equation*}
        \frac{1}{2a_j}\sum_{n = 1}^{2a_j}\PP\left(\tilde\gamma_1^\infty\cdots \tilde\gamma_n^\infty \cdot \ZZ^d \in K\right) \le 1 - \frac{\eps_p}{4},
    \end{equation*}
    which implies \eqref{eq:lim-transcendental}. 

    The rest of the proof is devoted to showing \eqref{eq.thm.1.3.goal}. By Markov's inequality and additivity of the expectation, we have: 
    \begin{equation*}
        \forall n \ge 1,\, \PP\left(\sum_{k = 1}^n H(g_k) \ge 4Mn/\eps_p \right) \le \frac{\eps_p}{4} 
 \end{equation*}
    Moreover, by Lemma \ref{lem:ex-lp}, we have:
    \begin{equation*}
        \forall n \ge 1,\; \PP\left(2\max_{k \le n} u_k^{-p} - \sum_{k = 1}^n u_k^{-p} \ge (2n)^p\right) \ge \eps_p.
    \end{equation*}
    Therefore, since $H$ is subadditive, $\tilde{\gamma}^j_n := (\gamma^j_n)^{s_n}$, $\gamma^j_n = h^j_n \gamma_n$, where $\gamma_n$ is as in \eqref{eq.def.gamma_n}, and $H(h^j_n) \le i_j \log(2)$ for all $j \ge 1$, we have:
    \begin{equation*}
        \forall n \ge 1,\; \PP\left(2\max_{k \le n} u_k^{-p} - \sum_{k = 1}^n H(\tilde\gamma_k^j) \ge (2n)^p - n \left(\frac{4M}{\eps_p} + i_j\log(2)\right)\right) \ge \frac{3\eps_p}{4}.
    \end{equation*}
    Therefore, by Lemma \ref{lem:systole} and using again $H(h^j_n) \le i_j \log(2)$, we deduce:
    \begin{equation*}
        \forall n \ge 1,\; \PP\left(-\log\delta\left(\tilde\gamma_1^j\cdots \tilde\gamma_n^j \cdot \ZZ^d\right) \ge (2n)^p - n \left(\frac{4M}{\eps_p} + 2i_j\log(2)\right) \right) \ge \frac{3\eps_p}{4}.
    \end{equation*}
    Therefore, by definition of $a_{j}$, we have:
    \begin{equation} \label{eq:bound-on-approx-systole}
        \forall n \ge a_j,\; \PP\left(-\log\delta\left(\tilde\gamma_1^j\cdots \tilde\gamma_n^j \cdot \ZZ^d\right) \ge (2^p - 1) n^p \right) \ge \frac{3\eps_p}{4}.
    \end{equation}

    By concavity of $x \mapsto x^{p'}$, we have for all $n \ge 1$:
    \begin{equation*}
        \int \log(\|\gamma\|\|\gamma^{-1}\|)^{p'} d\mu_j^{*n}(\gamma) \le n \int \log(\|\gamma\|\|\gamma^{-1}\|)^{p'} d\mu_j(\gamma).
    \end{equation*}
    Moreover, $\lambda_{j}$ is supported on $[0,1]$ so $\log(\|h^j_k\|\|(h_k^j)^{-1}\|) \le 1$ almost surely and for all $k$ and $j$. 
    Therefore $\int \log(\|\gamma\|\|\gamma^{-1}\|)^{p'} d\mu_j(\gamma) \le 1+ M'$. 
    Hence by Markov's inequality:
    \begin{equation*}
        \forall n \ge 1,\, \PP\left(\prod_{k = 1}^n \|\gamma^j_k\|\|(\gamma^j_k)^{-1}\| \ge \exp\left(\left(\frac{4n (M' + 1)}{\eps_p}\right)^{1/p'}\right)\right) \le \frac{\eps_p}{4}.
    \end{equation*}
    Therefore, by construction of $i_{j+1}$, for all $a_j \le n \le 2 a_j$, we have:
    \begin{equation*}
        \PP\left(\left(1 + 2^{-i_{j+1}}n (1+2^{-i_{j+1}})^n \prod_{k = 1}^n \|\gamma^j_k\|\|(\gamma^j_k)^{-1}\|\right) \ge 2\right) \le \frac{\eps_p}{4}.
    \end{equation*}
    Moreover, for all $x$, we have $\|N_x - I_d\| = |x|$ and $\|l_\infty(k_n) - l_j(k_n)\| \le 2^{-i_{j+1}+1}$. Therefore, we have $\frac{\|\tilde{\gamma}^\infty_k - \tilde{\gamma}^j_k\|}{\|\tilde{\gamma}^j_k\|} \le 2^{-i_{j+1}+1}$
    Hence, by Lemma \ref{lem:perturb-conclusion}, for all $n \ge 1$, we have:
    \begin{equation*}
        \frac{\delta\left(\tilde{\gamma}^\infty_1\cdots \tilde{\gamma}^\infty_n \cdot \ZZ^d\right)}{\delta\left(\gamma^j_1 \cdots \gamma^j_n \cdot \ZZ^d\right)} \le \left(1 + 2^{-i_{j+1}+1}n (1+2^{-i_{j+1}+1})^n \prod_{k = 1}^n \|\tilde\gamma^j_k\|\|(\tilde\gamma_k^j)^{-1}\|\right).
    \end{equation*}
    Hence, for all $a_j \le n \le 2 a_j$, we have:
    \begin{equation*}
        \PP\left({\delta\left(\tilde{\gamma}^\infty_1\cdots \tilde{\gamma}^\infty_n \cdot \ZZ^d\right)} \ge 2 {\delta\left(\gamma^j_1 \cdots \gamma^j_n \cdot \ZZ^d\right)}\right) \le \frac{\eps_p}{4}.
    \end{equation*}
    Then by \eqref{eq:bound-on-approx-systole}, we have:
    \begin{equation*}
        \forall a_j \le n \le 2a_j,\; \PP\left(-\log\delta\left(\tilde\gamma_1^\infty\cdots \tilde\gamma_n^\infty \cdot \ZZ^d\right) \ge (2^p - 1) n^p - \log(2) \right) \ge \frac{3\eps_p}{4} - \frac{\eps_p}{4},
    \end{equation*}
    which implies \eqref{eq.thm.1.3.goal} and hence finishes the proof.
\end{proof}

\subsection{Proof of Theorem \ref{th:non-alg-escape}}

We start by introducing the notion of a law with simple record. It will be used to construct the tail of the measure $\mu$ in Theorem \ref{th:non-alg-escape}.

\begin{Def}
    Let $\eta$ be a probability distribution on $\ZZ$ and let $(x_n) \sim \eta^{\otimes\NN}$. We say that $\eta$ has simple record if we have almost-surely:
    \begin{equation*}
        \#\left\{ k \le n \,\middle|\, x_k = \max_{i \le n} x_i\right\} \rightarrow 1 \;\;\; \text{as} \;\;\; n \to \infty.
    \end{equation*}
    In other words, with probability $1$, the maximum $\max_{i \le n} x_i$ is reached exactly once for all but finitely many $n$.
\end{Def}

The next lemma gives a construction of a law with simple record.

\begin{Lem}\label{lem:simple-record}
    Let $\eta$ be the measure on $\NN_{\ge 0}$ such that $\eta\{l\} = l^{-1/3} -(l+1)^{-1/3}$ \ie the push-forward of the uniform measure by $u \mapsto \lfloor u^{-3}\rfloor$.
    Then $\eta$ has simple record.
    Moreover, for $\eta^{\otimes\NN}$-a.e.\ $(x_i)$, we have $\max_{i \le n} x_i \ge n^2$ for all but finitely many $n$.
\end{Lem}

\begin{proof}
    Let $(x_i) \sim \eta^{\otimes\NN}$. 
    First we claim that $\max_{i \le n} x_i \ge n^{2}$ for all but finitely many $n$ almost surely.
    Indeed, for all $n \ge 1$, we have:
    \begin{align*}
        \PP\left(\max_{i \le n} x_i < n^{2}\right)  = \eta\{1, \dots, n^2-1\}^n  = \left(1- n^{-2/3}\right)^n \le \exp\left(-n^{1/3}\right).
    \end{align*} 
    Moreover $\sum_{n\ge 1} \exp\left(-n^{1/3}\right) < +\infty$ so by Borel--Cantelli Lemma, with probability $1$ there exist only finitely many indices $n$ such that $\max_{i \le n} x_i < n^{2}$.
    Now we claim that for all $1 \le i \neq j \le n$ and for all $l \in \NN_{\ge 0}$, we have:
    \begin{equation}\label{eq:i-neq-j}
        \PP\left(x_j = l\,\middle|\,  \max_{k \le n} x_k = x_i = l\right) = \frac{\eta\{l\}}{\eta\{0, \dots, l\}}.
    \end{equation} 
    Indeed, for all $l$, the conditional distribution of $(x_1, \dots, x_n)$ knowing that $\max_{k \le n} x_k \le l$ is $\frac{\mathds{1}_{\{1, \dots l\}^n}\eta^{\otimes n}}{\eta^{\otimes n}(\{0, \dots, l\}^n)} = \left(\frac{\mathds{1}_{\{0, \dots l\}}\eta}{\eta\{0, \dots, l\}}\right)^{\otimes n}$. 
    Therefore, for all $1 \le i \le n$, the conditional distribution of $(x_j)_{i \neq j \le n}$ knowing that $\max_{k \le n} x_k \le l$ and the value of $x_i$ is $\left(\frac{\mathds{1}_{\{0, \dots l\}}\eta}{\eta\{0, \dots, l\}}\right)^{\otimes n-1}$, which proves \eqref{eq:i-neq-j}.
    Note moreover that for all $l \in\NN$, we have $\eta\{l\} \le l^{-4/3}$ and $l \mapsto \eta\{l\}$ is a decreasing function on $\NN_{\ge 0}$.
    Moreover, for all $l \ge 7$, we have $\eta\{0, \dots, l\} \ge 1/2$.
    By \eqref{eq:i-neq-j}, we have for all $n \ge 3$, for all $l \ge n^2$ and for all $i \neq j \le n$:
    \begin{equation*}
        \PP\left(x_j = l\,\middle|\,  \max_{k \le n} x_k = x_i = l\right) \le 2n^{-8/3}.
    \end{equation*}
    Hence, for all $n \ge 3$, we have:
    \begin{equation*}
        \PP\left( \exists i \neq j \le n,\, x_i = x_j = \max_{k \le n} x_k \,\middle|\, \max_{k \le n} x_k \ge n^2 \right) \le (n-1) 2{n^{-8/3}} \le 2{n^{-5/3}}.
    \end{equation*}
    Hence, for all $n \ge 3$, we have:
    \begin{equation*}
        \PP\left( \exists i \neq j \le n,\, x_i = x_j = \max_{k \le n} x_k \right) \le 2{n^{-5/3}} + \exp\left(-n^{1/3}\right).
    \end{equation*}
    Moreover $\sum_{n\ge 3} 2{n^{-5/3}} < + \infty$ so by Borel Cantelli's Lemma, with probability $1$, the set of indices $n$ such that $\max_{i \le n} x_i$ is reached at least twice is almost surely finite, which concludes the proof.
\end{proof}

We are now ready to give the proof of Theorem \ref{th:non-alg-escape}. For the proof, we again start with a measure $\kappa_0$ that has full support and such that $H_*\kappa_0$ has finite first moment and then take its barycenter with a measure $\nu$ that is the push-forward of a measure with simple records by a well chosen function and convolute this with a measure $N_*{\lambda_{\infty}}$, which has no atoms.
However, we have to construct both $\nu$ and $\lambda_{\infty}$ at the same time.

\begin{proof}[Proof of Theorem \ref{th:non-alg-escape}]
    Let $\eta$ be the measure on $\NN_{\ge 0}$ such that $\eta\{l\} = l^{-1/3} -(l+1)^{-1/3}$.
    Let $\kappa$ be a probability measure on $\SL_d(\QQ)$ that has full support and such that $H_*\kappa$ has finite first moment and write $M := \int H d\kappa$.
    We define two increasing sequences of integers $(i_j)_{j \ge 1}$ and $(l_j)_{j \ge 1}$ by induction.
    Let $l_1 = i_1 = 1$.
    Let $j \ge 1$ and assume that we have constructed $i_j$ and $l_j$.
    Assume also that $l_j \le i_j$.
    Let $l_{j+1} := \lceil j (i_j\log(2) + l_j \log(2) + 2M)\rceil$ and let $i_{j+1}$ be the smallest integer such that $i_{j+1} \ge \max\{l_j, i_j + 1\}$ and:
    \begin{equation*}
        2^{-i_{j+1}+1} j (1 + 2^{-i_{j+1}+1})^{j} 2\exp((2M + l_{j+1} \log(2) + 1) j) \le 1.
    \end{equation*}
    Let $(j_n) \sim \eta^{\otimes \NN}$ and let $(g_n) \sim \kappa^{\otimes\NN}$ be independent.
    Let $\left(\epsilon^n_k\right)_{n \ge 0, k \ge 1}\in\{0,1\}^{\NN^2}$ and $(s_n)\in\{-1,1\}^\NN$ be fixed.
    For all $n \ge 0$, we write $t_n := \sum_{k \geq 1} \epsilon_n^k 2^{-i_k}$.
    For all $n \ge 1$, we write $\gamma_n := g_n$ when $j_n = 0$ and $\gamma_n = \diag(2^{l_{j_n}}, 2^{-l_{j_n}}, 1, \dots, 1)$ otherwise, and write $\tilde{\gamma}_n = (N_{t_n} \gamma_n)^{s_n}$.
    For all $j,n \in \NN$, we write $t_n^j := \sum_{k=1}^{j} \epsilon_n^k 2^{-i_k}$ and $h_n^j := N_{t_n^j}$ and $\tilde{\gamma}_n^j := (h_n^j \gamma_n)^{s_n}$. Note that for all $j,n$, we have $|t_n - t_n^j| \le 2^{-i_{j+1}+1}$ and $\frac{\|\tilde{\gamma}_n^j - \tilde{\gamma}_n\|}{\|\tilde{\gamma}_n\|} \le 2^{-i_{j+1}+1}$.

    Write $n_0$ for the random smallest integer such that for all $n \ge n_0$, we have $n^2 \le \max_{k' \le n} j_{k'}$ and $\#\left\{k \le n\,\middle|\, j_k = \max_{k' \le n} j_{k'}\right\} = 1$ and $\sum_{k = 1}^{n} H(g_k) \le 2M n$.
    Such a $n_0$ exists with probability $1$ by Lemma \ref{lem:simple-record} and by the law of large numbers.
    Let $(j_n)$ and $(g_n)$ be fixed such that $n_0$ is finite and let $n \ge \max\{n_0, 2\}$.
    Let $k_n$ be the unique integer such that $j_{k_n} = \max_{k' \le n} j_{k'}$ and let $j := j_{k_n} - 1$.
    Then we have $\sum_{k \le n}H(g_k) \le 2M n$, therefore $\sum_{k \le n} H(\gamma_k) \le l_{j+1} + n l_{j} \log(2) + 2Mn$ and $\sum_{k \le n} H(h_k^{j}) \le n i_{j} \log(2)$. 
    So by Lemma \ref{lem:systole}, we have:
    \begin{multline*}
        -\log\delta\left(\tilde{\gamma}^j_1 \cdots \tilde{\gamma}^j_n h^j_0 \ZZ^d\right) \ge l_{j + 1} - n(l_j \log(2) + i_j\log(2) + 2 M) \\ \ge \left(1 - \frac{n}{n^2-1}\right) l_{j+1} \ge \frac{l_{j+1}}{3}.
    \end{multline*}
    Moreover, for all $k \le n$, we have $\|h_k^j\|, \|(h_k^j)^{-1}\| \le 2$ and $\|\gamma_k\|\le \max\{2M, l_{j+1}\}$. 
    By Lemma \ref{lem:perturb-conclusion}, we have:
    \begin{equation*}
        \frac{\delta\left(\tilde{\gamma}_1 \cdots \tilde{\gamma}_nh_0 \ZZ^d\right)}{\delta\left(\tilde{\gamma}^j_1 \cdots \tilde{\gamma}^j_n h_0^j \ZZ^d\right)} \le 1+ 2^{-i_{j+1}+1} j (1 + 2^{-i_{j+1}+1})^{j} 2\exp((2M + l_{j+1} \log(2) + 1)n) \le 2.
    \end{equation*}
    From that, we conclude that $\delta\left(\tilde{\gamma}_1 \cdots \tilde{\gamma}_n h_0 \ZZ^d\right) \to 0$.

    Now let $(j_n) \sim \eta^{\otimes \NN}$, let $(g_n) \sim \kappa^{\otimes\NN}$ and
    let $\left(\epsilon^n_k\right)_{n \ge 1, k \ge 1}\in\{0,1\}^{\NN^2}$ and $(s_n)\in\{-1,1\}^\NN$ be uniformly distributed and assume that all these random variables are globally independent.
    Let $S := \{N_t\,|\, t = \sum_{k = 1}^\infty \epsilon_k 2^{-i_k}, (\epsilon_k)\in\{0,1\}^\NN\}$.
    The set $S$ is in bijection with $\{0,1\}^\NN$ so it is uncountable.
    Then $(\tilde{\gamma}_n)_n$ is i.i.d., write $\mu$ for the distribution of $\tilde\gamma_1$.
    Then $\mu$ is symmetric has full support and no atoms and $\delta\left(\tilde{\gamma}_1 \cdots \tilde{\gamma}_n h_0 \ZZ^d\right) \to 0$ for all $h_0 \in S$ with probability $1$, as desired.
\end{proof}

We will note that on the above example we have no control over the moments of $\mu$ nor on its regularity.

\bibliographystyle{alpha}

\bibliography{biblio.bib}

\end{document}